\documentclass[twoside,reqno,11pt]{amsart}
\usepackage{latexsym}%%triangles:\lhd\rhd\unlhd\unrhd%
\newcommand{\qdn}{\hspace*{-1.5mm}}
\newcommand{\qqdn}{\hspace*{-2.5mm}}
\newcommand{\xqdn}{\hspace*{-5.0mm}}
\newcommand{\xxqdn}{\hspace*{-10mm}}

\newcommand{\fns}{\footnotesize}

\newcommand{\ffnk}[4]{\left[\qdn\ba{#1}#3\\[0.7mm]#4\ea{\!\Big|\:#2}\right]}

\newcommand{\nnm}{\nonumber}
\newcommand{\be}{\begin{equation}}
\newcommand{\ee}{\end{equation}}
\newcommand{\ba}{\begin{array}}
\newcommand{\ea}{\end{array}}
\newcommand{\bmn}{\begin{eqnarray}}
\newcommand{\emn}{\end{eqnarray}}
\newcommand{\bnm}{\begin{eqnarray*}}
\newcommand{\enm}{\end{eqnarray*}}
\newcommand{\bln}{\begin{subequations}}
\newcommand{\eln}{\end{subequations}}

\newtheorem{thm}{Theorem}%[section]

\newtheorem{exam}{Example}

\newtheorem{entry}{Entry}%%%%%%%%%%%%%%%%

\newcommand{\bbtm}[4]{\bibitem{kn:#1}{#2,}~{#3,}~{#4.}}
\newcommand{\cito}[1]{\cite{kn:#1}}
\newcommand{\citu}[2]{\cite[#2]{kn:#1}}

\begin{document} %%%%%%%%%% This paper is published in %%%%%%%
{\fns% \today\hfill\copyright%% Printed in China} %%%%%%%%%%%%%%%
%%%%%%%%%%%%%%%%%%%%%%%%%%%%%%%%%%%%%%%%%%%%%%%%%%%%%%%%%%%%%%
\title{$\pi$-formulas with free parameters}
\author{$^a$Chuanan Wei, $^b$Dianxuan Gong, $^c$ Jianbo Li}
\dedicatory{
$^A$Department of Information Technology\\
  Hainan Medical College,  Haikou 571101, China\\
         $^B$College of Sciences\\
             Hebei Polytechnic University, Tangshan 063009, China \\
   $^C$Institute of Mathematical Sciences\\
 Xuzhou Normal University, Xuzhou 221116, China}
\thanks{\emph{Email addresses}:
 weichuanan@yahoo.com.cn (C. Wei), gongdianxuan@yahoo.com.cn (D. Gong),
 ljianb66@gmail.com (J. Li)}
\address{ }


\begin{abstract}
In terms of the hypergeometric method, we give the extensions of two
known series for $\pi$. Further, other twenty-nine summation
formulas for $\pi$, $\pi^2$ and $1/\pi$ with free parameters are
also derived in the same way.
\end{abstract}

%%%%%%%%%%%%%%%%%%%%%%%%%%%%%%%%%%%%%%%%%%%%%%%%%%%%%%%%%%%%%%%%%%%
\maketitle\thispagestyle{empty}%%%%%%%%%%%%%%%%%%%%%%%%%%%%%%%%%%%%
\markboth{Chuanan Wei, Dianxuan Gong}%%%%%%%%%%%%%%%%%%%%%%%%%%%%
         {$\pi$-formulas with free parameters}

%%%%%%%%%%%%%%%%%%%%%%%%%%%%%%%%%%%%%%%%%%%%%%%%%%%%%%%%%%%%%%%%%%%
%%%%%%%%%%%%%%%%%%%%%%%%%%%%%%%%%%%%%%%%%%%%%%%%%%%%%%%%%%%%%%%%%%%
%%%%%%%%%%%%%%%%%%%%%%%%%%%%%%%%%%%%%%%%%%%%%%%%%%%%%%%%%%%%%%%%%%%
\section{Introduction}

 For a complex number $x$ and an integer $n$, define the shifted
factorial by
 \bnm
 (x)_n=
\begin{cases}
\prod_{i=0}^{n-1}(x+i),&\quad\text{when}\quad n>0;\\
1,&\quad\text{when}\quad n=0;\\
\frac{1}{\prod_{j=n}^{-1}(x+j)},&\quad\text{when}\quad n<0.
\end{cases}
 \enm
Define $\Gamma$-function by Euler's integral:
\[\Gamma(x)=\int_{0}^{\infty}t^{x-1}e^{-t}dt\:\:\text{with}\:\:Re(x)>0.\]
Then we have the following two relations:
 \bnm
 \Gamma(x+n)=\Gamma(x)(x)_n,\quad\Gamma(x)\Gamma(1-x)
  =\frac{\pi}{\sin(\pi x)},
  \enm
 which will frequently be used without indication in this paper.

Following  Bailey~\cito{bailey}, define the hypergeometric series by
\[_{1+r}F_s\ffnk{cccc}{z}{a_0,&a_1,&\cdots,&a_r}
{&b_1,&\cdots,&b_s} \:=\:\sum_{k=0}^\infty
\frac{(a_{0})_{k}(a_{1})_{k}\cdots(a_{r})_{k}}
 {k!(b_{1})_{k}\cdots(b_{s})_{k}}z^k.\]
Then a simple $_2F_1$-series identity (cf. \citu{weisstein}{Eq.
(26)}) can be stated as
 \bnm
_2F_1\ffnk{ccccc}{x^2}{1,&1}
 {&\frac{3}{2}}=\frac{\arcsin(x)}
 {x\sqrt{1-x^2}}\quad\text{where}\quad |x|<1.
 \enm
Two beautiful series for $\pi$ (cf. \citu{weisstein}{Eq. (23) and
(27)}) implied by it read as
  \bmn
\label{known-a}
 &&\xxqdn\sum_{k=0}^{\infty}\frac{k!}{(2k+1)!!}=\frac{\pi}{2},\\
\label{known-b}
 &&\xxqdn\sum_{k=0}^{\infty}\frac{(k!)^2}{(2k+1)!}=\frac{2\pi}{3\sqrt{3}}.
 \emn

Recall the $_7F_6$-series identity due to Chu \citu{chu-a}{Eq.
(5.1e)} and Dougall's $_5F_4$-series identity (cf. \citu{bailey}{p.
27}):
 \bmn
 &&\xqdn_7F_6\ffnk{cccccccc}{1}{a-\frac{1}{2},&\frac{2a+2}{3},&2b-1,&2c-1,&2+2a-2b-2c,&a+s,&-s}
 {&\frac{2a-1}{3},&1+a-b,&1+a-c,&b+c-\frac{1}{2},&2a+2s,&-2s}\nnm\\
 &&\xqdn\:\,=\:\,\frac{(\frac{1}{2}+a)_s(b)_s(c)_s(a-b-c+\frac{3}{2})_s}
 {(\frac{1}{2})_s(1+a-b)_s(1+a-c)_s(b+c-\frac{1}{2})_s},
 \label{hypergeotric-b}\\
&&\xqdn_5F_4\ffnk{cccccccc}{1}{a,&1+\frac{a}{2},&b,&c,&d}
 {&\frac{a}{2},&1+a-b,&1+a-c,&1+a-d}\nnm\\
 &&\xqdn\:\,=\:\,\frac{\Gamma(1+a-b)\Gamma(1+a-c)\Gamma(1+a-d)\Gamma(1+a-b-c-d)}
 {\Gamma(1+a)\Gamma(1+a-b-c)\Gamma(1+a-b-d)\Gamma(1+a-c-d)}
 \label{hypergeotric-a}
 \emn
provided that $Re(1+a-b-c-d)>0$.

Recently, Chu \cite{kn:chu-b,kn:chu-c} and Liu
\cite{kn:liu-a,kn:liu-b} deduced many surprising $\pi$-formulas from
several known hypergeometric series identities. Thereinto, Chu
\cito{chu-a} showed that \eqref{hypergeotric-a} implies the
Ramanujan-type series for $1/\pi$ with three free parameters:
 \bmn \label{chu}\qquad
 \frac{2}{\pi}=\frac{(\frac{1}{2})_{m-n-p}}{(\frac{1}{2})_{n}(\frac{1}{2})_{p}}
\sum_{k=0}^{\infty}(-1)^k
 \frac{(\frac{1}{2})_{k+m}(\frac{1}{2})_{k+n}(\frac{1}{2})_{k+p}}
 {k!(k+m-n)!(k+m-p)!}(4k+2m+1)
 \emn
where $m,n,p\in \mathbb{Z}$ with $\min\{m-n,m-p,m-2n-2p\}\geq0$ and
Liu \cito{liu-b} showed that \eqref{hypergeotric-a} implies the
Ramanujan-type series for $1/\pi$ with four free parameters:
 \bmn
\quad\frac{\sqrt{3}}{3\pi}&&\xqdn=\:\frac{(\frac{2}{3})_{m-n-p}(\frac{1}{3})_{m-n-q}(\frac{1}{2})_{m-p-q}}
 {(m-n-p-q-1)!(\frac{1}{2})_{n}(\frac{1}{3})_{p}(\frac{2}{3})_{q}}
\nnm \\\label{liu} &&\xqdn\times\:\: \sum_{k=0}^{\infty}
 \frac{(\frac{1}{2})_{k+m}(\frac{1}{2})_{k+n}(\frac{1}{3})_{k+p}(\frac{2}{3})_{k+q}}
 {k!(k+m-n)!(\frac{7}{6})_{k+m-p}(\frac{5}{6})_{k+m-q}}(1+2m+4k)
 \emn
where $m,n,p,q\in \mathbb{Z}$ with $\min\{m-n,m-n-p-q-1\}\geq0$.
More $\pi$-formulas can be found in the papers
\cito{adanmchik}-\cito{bailey-b}, \cito{baruch-a}-\cito{chan-d},
\cito{glaisher}, \cito{goura}-\cito{levrie} and
\cito{ramanujan}-\cito{zudilin-b}.

Inspired by these work just mentioned, we shall explore further the
relations of $\pi$-formulas and hypergeometric series. The structure
of the paper is arranged as follows. Seven $\pi$-formulas with free
parameters including the extensions of \eqref{known-a} and
\eqref{known-b} will be derived from \eqref{hypergeotric-b} in
section 2. Twenty-four $\pi$-formulas with free parameters  other
than \eqref{chu} and \eqref{liu} will be deduced from
\eqref{hypergeotric-a} in section 3.

%%%%%%%%%%%%%%%%%%%%%%%%%%%%%%%%%%%%%%%%%%%%%%%%%%%%%%%%%%%%%%%%%%%
%%%%%%%%%%%%%%%%%%%%%%%%%%%%%%%%%%%%%%%%%%%%%%%%%%%%%%%%%%%%%%%%%%%
\section{Summation formulas for $\pi$ and  $\pi^2$ with free parameters \\
implied by Chu's $_7F_6$-series identity}
%%%%%%%%%%%%%%%%%%%%%%%%%%%%%%%%%%%%%%%%%%%%%%%%%%%%%%%%%%%%%%%%%%%
%%%%%%%%%%%%%%%%%%%%%%%%%%%%%%%%%%%%%%%%%%%%%%%%%%%%%%%%%%%%%%%%%%%

Letting $s\to\infty$ for \eqref{hypergeotric-b}, we obtain the
following equation:
 \bmn
 &&\xqdn\qqdn_5F_4\ffnk{cccccccc}{\frac{1}{4}}{a-\frac{1}{2},&\frac{2a+2}{3},&2b-1,&2c-1,&2+2a-2b-2c}
 {&\frac{2a-1}{3},&1+a-b,&1+a-c,&b+c-\frac{1}{2}}\nnm\\ \label{hypergeotric-c}
 &&\xqdn\qqdn\:\,=\:\,
 \frac{\Gamma(\frac{1}{2})\Gamma(1+a-b)\Gamma(1+a-c)\Gamma(b+c-\frac{1}{2})}
 {\Gamma(\frac{1}{2}+a)\Gamma(b)\Gamma(c)\Gamma(a-b-c+\frac{3}{2})}.
 \emn
Subsequently, one summation formula for $\pi$ with two free
parameters, five summation formulas for $\pi$ with three free
parameters and one summation formula for $\pi^2$ with three free
parameters will respectively be derived from \eqref{hypergeotric-c}.

Choosing $b=1+m$, $c=1+n$ in \eqref{hypergeotric-c} and then letting
$a\to\infty$, we achieve the extension of \eqref{known-a}.

\begin{thm}\label{thm-one} For $m,n\in \mathbb{N}_0$,
 there holds the summation formula for $\pi$ with two free
parameters:
\bnm
 \frac{\pi}{2^{m+n+1}}=\frac{m!n!}{(2m)!(2n)!}\sum_{k=0}^{\infty}\frac{(k+2m)!(k+2n)!}{k!(2k+2m+2n+1)!!}.
 \enm
\end{thm}

When $m=n=0$, Theorem \ref{thm-one} reduces to \eqref{known-a}
exactly. Other two examples of the same type are displayed as
follows.

\begin{exam}[$m=1$ and $n=0$ in Theorem \ref{thm-one}]
 \bnm\:\:
 \frac{\pi}{2}=\sum_{k=1}^{\infty}\frac{(k+1)!}{(2k+1)!!}.
 \enm
\end{exam}

\begin{exam}[$m=2$ and $n=0$ in Theorem \ref{thm-one}]
 \bnm\quad
 \frac{3\pi}{2}=\sum_{k=2}^{\infty}\frac{(k+2)!}{(2k+1)!!}.
 \enm
\end{exam}

Making $a=\frac{1}{2}+m$, $b=\frac{1}{4}+n$ and $c=\frac{3}{4}+p$ in
\eqref{hypergeotric-c}, we attain the equation.

\begin{thm}\label{thm-two} For $m,n,p\in \mathbb{Z}$ with $\min\{m,m-n-p\}\geq0$,
 there holds the summation formula for $\pi$ with three free
parameters:
 \bnm \quad
 \frac{\pi}{2}&&\xqdn=\:\frac{1}{(-\frac{1}{4})_{n}(\frac{1}{4})_{p}(\frac{1}{2})_{m-n-p}}
\\[0.5mm]&&\xqdn\times\:\: \sum_{k=0}^{\infty}
 \frac{(-\frac{1}{2})_{k+2n}(\frac{1}{2})_{k+2p}(k+m)!(k+2m-2n-2p)!}
 {(\frac{5}{4})_{k+m-n}(\frac{3}{4})_{k+m-p}(\frac{1}{2})_{k+n+p}\,k!}\frac{3k+2m}{k+m}\frac{1}{4^{k+m}}.
 \enm
\end{thm}

Two examples from Theorem \ref{thm-two} are laid out as follows.

\begin{exam}[$m=n=1$ and $p=0$ in Theorem \ref{thm-two}]
 \bnm\:\:
 \frac{\pi}{4}=\sum_{k=0}^{\infty}\frac{(2k)!}{(4k+3)!!}(3k+2).
 \enm
\end{exam}

\begin{exam}[$m=2$ and $n=p=1$ in Theorem \ref{thm-two}]
 \bnm\:\:
 \frac{\pi}{4}=\sum_{k=1}^{\infty}\frac{(2k)!}{(4k+1)!!}(3k+1).
 \enm
\end{exam}

A beautiful result should be mentioned. Fixing $a=-\frac{3}{2}$ in
the identity due to Chu \citu{chu-a}{Equation (5.3 g)}:
 \bnm
&&_6F_5\ffnk{cccccccc}{\frac{32}{27}}{7+4a,&-2-2a,&-\frac{5}{2}-2a,&\frac{2-4a}{5},&-a+s,&-s}
 {&\frac{4}{3},&\frac{5}{3},&-\frac{3+4a}{5},&-2-4s,&-2-4a+4s}\nnm\\
 &&\:\,=\:\,\frac{(-\frac{1}{2}-a)_s(-\frac{1}{4}-a)_s(\frac{1}{4}-a)_s(\frac{5}{2}+a)_s}
 {(\frac{3}{4})_s(\frac{5}{4})_s(\frac{3}{2})_s(-\frac{3}{2}-2a)_s}
 \enm
and then letting $s\to\infty$, we obtain the surprising series for
$\pi$:
 \bnm\quad
 \frac{\pi}{2}=\sum_{k=0}^{\infty}\frac{k!(2k)!}{(3k+2)!}\frac{5k+3}{2^k}.
 \enm

Setting $a=\frac{1}{2}+m$, $b=\frac{1}{6}+n$ and $c=\frac{5}{6}+p$
in \eqref{hypergeotric-c}, we get the extension of \eqref{known-b}.

\begin{thm}\label{thm-three} For $m,n,p\in \mathbb{Z}$ with $\min\{m,m-n-p\}\geq0$,
 there holds the summation formula for $\pi$ with three free
parameters:
 \bnm \quad
 \frac{2\pi}{3\sqrt{3}}&&\xqdn=\:\frac{1}{(-\frac{1}{3})_{n}(\frac{1}{3})_{p}(\frac{1}{2})_{m-n-p}}
\\[0.5mm]&&\xqdn\times\:\: \sum_{k=0}^{\infty}
 \frac{(-\frac{2}{3})_{k+2n}(\frac{2}{3})_{k+2p}(k+m)!(k+2m-2n-2p)!}
 {(\frac{4}{3})_{k+m-n}(\frac{2}{3})_{k+m-p}(\frac{1}{2})_{k+n+p}\,k!}\frac{3k+2m}{k+m}\frac{1}{4^{k+m}}.
 \enm
\end{thm}

When $m=n=1$ and $p=0$, Theorem \ref{thm-three} reduces to
\eqref{known-b} exactly. Other two examples of the same type are
displayed as follows.

\begin{exam}[$m=2$ and $n=p=1$ in Theorem \ref{thm-three}]
 \bnm
 \frac{2\pi}{3\sqrt{3}}=\sum_{k=1}^{\infty}\frac{(k!)^2}{(2k+1)!}(3k+2).
 \enm
\end{exam}

\begin{exam}[$m=3$, $n=2$ and $p=1$ in Theorem \ref{thm-three}]
 \bnm\:\:
 \frac{4\pi}{9\sqrt{3}}=\sum_{k=2}^{\infty}\frac{(k!)^2}{(2k+1)!}(3k+1)k.
 \enm
\end{exam}

Taking $a=\frac{1}{2}+m$, $b=\frac{1}{3}+n$ and $c=\frac{2}{3}+p$ in
\eqref{hypergeotric-c}, we gain the equation.

\begin{thm}\label{thm-four} For $m,n,p\in \mathbb{Z}$ with $\min\{m,m-n-p\}\geq0$,
 there holds the summation formula for $\pi$ with three free
parameters:
 \bnm \quad
 \frac{\pi}{\sqrt{3}}&&\xqdn=\:\frac{1}{(-\frac{1}{6})_{n}(\frac{1}{6})_{p}(\frac{1}{2})_{m-n-p}}
\\[0.5mm]&&\xqdn\times\:\: \sum_{k=0}^{\infty}
 \frac{(-\frac{1}{3})_{k+2n}(\frac{1}{3})_{k+2p}(k+m)!(k+2m-2n-2p)!}
 {(\frac{7}{6})_{k+m-n}(\frac{5}{6})_{k+m-p}(\frac{1}{2})_{k+n+p}\,k!}\frac{3k+2m}{k+m}\frac{1}{4^{k+m}}.
 \enm
\end{thm}

Two examples from Theorem \ref{thm-four} are laid out as follows.

\begin{exam}[$m=n=1$ and $p=0$ in Theorem \ref{thm-four}]
 \bnm\:\:
 \frac{5\pi}{4\sqrt{3}}=\sum_{k=0}^{\infty}\frac{(1)_k(\frac{1}{3})_k(\frac{5}{3})_k}
 {(\frac{3}{2})_k(\frac{7}{6})_k(\frac{11}{6})_k}\frac{3k+2}{4^k}.
 \enm
\end{exam}

\begin{exam}[$m=2$ and $n=p=1$ in Theorem \ref{thm-four}]
  \bnm\:\:
 \frac{\pi}{2\sqrt{3}}=\sum_{k=1}^{\infty}\frac{(1)_k(\frac{2}{3})_k(\frac{4}{3})_k}
 {(\frac{3}{2})_k(\frac{5}{6})_k(\frac{7}{6})_k}\frac{3k+1}{4^k}.
 \enm
\end{exam}

Choosing $a=\frac{1}{2}+m$, $b=\frac{5}{12}+n$ and
$c=\frac{7}{12}+p$ in \eqref{hypergeotric-c}, we achieve the
equation.

\begin{thm}\label{thm-five} For $m,n,p\in \mathbb{Z}$ with $\min\{m,m-n-p\}\geq0$,
 there holds the summation formula for $\pi$ with three free
parameters:
\bnm \quad
 \frac{\pi}{6(2-\sqrt{3})}&&\xqdn=\:\frac{1}{(-\frac{1}{12})_{n}(\frac{1}{12})_{p}(\frac{1}{2})_{m-n-p}}
\\[0.5mm]&&\xqdn\times\:\: \sum_{k=0}^{\infty}
 \frac{(-\frac{1}{6})_{k+2n}(\frac{1}{6})_{k+2p}(k+m)!(k+2m-2n-2p)!}
 {(\frac{13}{12})_{k+m-n}(\frac{11}{12})_{k+m-p}(\frac{1}{2})_{k+n+p}\,k!}\frac{3k+2m}{k+m}\frac{1}{4^{k+m}}.
 \enm
\end{thm}

Two examples from Theorem \ref{thm-five} are displayed as follows.

\begin{exam}[$m=n=1$ and $p=0$ in Theorem \ref{thm-five}]
 \bnm\:\:
 \frac{11\pi}{60(2-\sqrt{3})}=\sum_{k=0}^{\infty}\frac{(1)_k(\frac{1}{6})_k(\frac{11}{6})_k}
 {(\frac{3}{2})_k(\frac{13}{12})_k(\frac{23}{12})_k}\frac{3k+2}{4^k}.
 \enm
\end{exam}

\begin{exam}[$m=2$ and $n=p=1$ in Theorem \ref{thm-five}]
  \bnm\:\:
 \frac{\pi}{12(2-\sqrt{3})}=\sum_{k=1}^{\infty}\frac{(1)_k(\frac{5}{6})_k(\frac{7}{6})_k}
 {(\frac{3}{2})_k(\frac{11}{12})_k(\frac{13}{12})_k}\frac{3k+1}{4^k}.
 \enm
\end{exam}

Making $a=\frac{1}{2}+m$, $b=\frac{1}{12}+n$ and $c=\frac{11}{12}+p$
in \eqref{hypergeotric-c}, we attain the equation.

\begin{thm}\label{thm-six} For $m,n,p\in \mathbb{Z}$ with $\min\{n,p,m-n-p\}\geq0$,
 there holds the summation formula for $\pi$ with three free
parameters:
 \bnm \quad
 \frac{5\pi}{6(2+\sqrt{3})}&&\xqdn=\:\frac{1}{(-\frac{5}{12})_{n}(\frac{5}{12})_{p}(\frac{1}{2})_{m-n-p}}
\\[0.5mm]&&\xqdn\times\:\: \sum_{k=0}^{\infty}
 \frac{(-\frac{5}{6})_{k+2n}(\frac{5}{6})_{k+2p}(k+m)!(k+2m-2n-2p)!}
 {(\frac{17}{12})_{k+m-n}(\frac{7}{12})_{k+m-p}(\frac{1}{2})_{k+n+p}\,k!}\frac{3k+2m}{k+m}\frac{1}{4^{k+m}}.
 \enm
\end{thm}

Two examples from Theorem \ref{thm-six} are laid out as follows.

\begin{exam}[$m=n=1$ and $p=0$ in Theorem \ref{thm-six}]
  \bnm\:\:
 \frac{35\pi}{12(2+\sqrt{3})}=\sum_{k=0}^{\infty}\frac{(1)_k(\frac{5}{6})_k(\frac{7}{6})_k}
 {(\frac{3}{2})_k(\frac{17}{12})_k(\frac{19}{12})_k}\frac{3k+2}{4^k}.
 \enm
\end{exam}

\begin{exam}[$m=2$ and $n=p=1$ in Theorem \ref{thm-five}]
  \bnm\:\:
 \frac{5\pi}{12(2+\sqrt{3})}=\sum_{k=1}^{\infty}\frac{(1)_k(\frac{1}{6})_k(\frac{11}{6})_k}
 {(\frac{3}{2})_k(\frac{7}{12})_k(\frac{17}{12})_k}\frac{3k+1}{4^k}.
 \enm
\end{exam}

Setting $a=\frac{3}{2}+m$, $b=1+n$ and $c=1+p$ in
\eqref{hypergeotric-c}, we get the equation.

\begin{thm}\label{thm-seven} For $m,n,p\in \mathbb{Z}$ with $\min\{n,p,m-n-p\}\geq0$,
 there holds the summation formula for $\pi^2$ with three free parameters:
\bnm \quad
 \pi^2&&\xqdn=\:\frac{1}{(\frac{1}{2})_{n}(\frac{1}{2})_{p}(\frac{1}{2})_{m-n-p}}
\\[0.5mm]&&\xqdn\times\:\: \sum_{k=0}^{\infty}
 \frac{(k+m)!(k+2n)!(k+2p)!(k+2m-2n-2p)!}
 {(\frac{3}{2})_{k+m-n}(\frac{3}{2})_{k+m-p}(\frac{3}{2})_{k+n+p}\,k!}\frac{3k+2m+2}{4^{k+m-1}}.
 \enm
\end{thm}

Two examples from Theorem \ref{thm-seven} are displayed as follows.

\begin{exam}[$m=n=p=0$ in Theorem \ref{thm-seven}]
 \bnm\:\:
 \frac{\pi^2}{4}=\sum_{k=0}^{\infty}\frac{(1)_k^3}
 {(\frac{3}{2})_k^3}\frac{3k+2}{4^{k}}.
 \enm
\end{exam}

\begin{exam}[$m=2$ and $n=p=1$ in Theorem \ref{thm-seven}]
  \bnm\:\:
 \frac{\pi^2}{12}=\sum_{k=2}^{\infty}\frac{(1)_k^3}
 {(\frac{1}{2})_k^2(\frac{3}{2})_k}\frac{k}{4^{k}}.
  \enm
\end{exam}
%%%%%%%%%%%%%%%%%%%%%%%%%%%%%%%%%%%%
%%%%%%%%%%%%%%%%%%%%%%%%%%%%%%%%%%%%%%%%%%%%%%%%%%%%%%%%%%%%%%%%%%%
%%%%%%%%%%%%%%%%%%%%%%%%%%%%%%%%%%%%%%%%%%%%%%%%%%%%%%%%%%%%%%%%%%%
\section{Ramanujan-type series for $1/\pi$ with free parameters\\implied by Dougall's $_5F_4$-series identity}
%%%%%%%%%%%%%%%%%%%%%%%%%%%%%%%%%%%%%%%%%%%%%%%%%%%%%%%%%%%%%%%%%%%
%%%%%%%%%%%%%%%%%%%%%%%%%%%%%%%%%%%%%%%%%%%%%%%%%%%%%%%%%%%%%%%%%%%
%%%%%%%%%%%%%%%%%%%%%%%%%%%%%%%%%%%%%%%%%%%%%%%%%%%%%%%%%%%%%%%%%%%

In this section, ten Ramanujan-type series for $1/\pi$ with three
free parameters and fourteen Ramanujan-type series for $1/\pi$ with
four free parameters other than \eqref{chu} and \eqref{liu} will
respectively be deduced from \eqref{hypergeotric-a}.

Taking $a=\frac{1}{6}+m$, $b=\frac{1}{6}+n$, $c=\frac{1}{6}+p$ in
\eqref{hypergeotric-a} and then letting $d\to\infty$, we gain the
equation.

\begin{thm}\label{thm-aa} For $m,n,p\in \mathbb{Z}$ with $\min\{m-n,m-p,m-2n-2p+1\}\geq0$,
 there holds the Ramanujan-type series for $1/\pi$ with three free parameters:
\bnm\qquad
 \frac{3}{\pi}=\frac{(\frac{5}{6})_{m-n-p}}{(\frac{1}{6})_{n}(\frac{1}{6})_{p}}
\sum_{k=0}^{\infty}(-1)^k
 \frac{(\frac{1}{6})_{k+m}(\frac{1}{6})_{k+n}(\frac{1}{6})_{k+p}}
 {k!(k+m-n)!(k+m-p)!}(1+6m+12k).
 \enm
\end{thm}

Choosing $a=\frac{5}{6}+m$, $b=\frac{5}{6}+n$, $c=\frac{5}{6}+p$ in
\eqref{hypergeotric-a} and then letting $d\to\infty$, we achieve the
equation.

\begin{thm}\label{thm-bb} For $m,n,p\in \mathbb{Z}$ with $\min\{m-n,m-p,m-2n-2p-1\}\geq0$,
 there holds the Ramanujan-type series for $1/\pi$ with three free parameters:
\bnm\qquad
 \frac{3}{\pi}=\frac{(\frac{1}{6})_{m-n-p}}{(\frac{5}{6})_{n}(\frac{5}{6})_{p}}
\sum_{k=0}^{\infty}(-1)^k
 \frac{(\frac{5}{6})_{k+m}(\frac{5}{6})_{k+n}(\frac{5}{6})_{k+p}}
 {k!(k+m-n)!(k+m-p)!}(5+6m+12k).
 \enm
\end{thm}

Making $a=\frac{1}{4}+m$, $b=\frac{1}{4}+n$, $c=\frac{1}{4}+p$ in
\eqref{hypergeotric-a} and then letting $d\to\infty$, we attain the
equation.

\begin{thm}\label{thm-cc} For $m,n,p\in \mathbb{Z}$ with $\min\{m-n,m-p,m-2n-2p+1\}\geq0$,
 there holds the Ramanujan-type series for $1/\pi$ with three free parameters:
\bnm\qquad
 \frac{2\sqrt{2}}{\pi}=\frac{(\frac{3}{4})_{m-n-p}}{(\frac{1}{4})_{n}(\frac{1}{4})_{p}}
\sum_{k=0}^{\infty}(-1)^k
 \frac{(\frac{1}{4})_{k+m}(\frac{1}{4})_{k+n}(\frac{1}{4})_{k+p}}
 {k!(k+m-n)!(k+m-p)!}(1+4m+8k).
 \enm
\end{thm}

Setting $a=\frac{3}{4}+m$, $b=\frac{3}{4}+n$, $c=\frac{3}{4}+p$ in
\eqref{hypergeotric-a} and then letting $d\to\infty$, we get the
equation.

\begin{thm}\label{thm-dd} For $m,n,p\in \mathbb{Z}$ with $\min\{m-n,m-p,m-2n-2p-1\}\geq0$,
 there holds the Ramanujan-type series for $1/\pi$ with three free parameters:
\bnm\qquad
 \frac{2\sqrt{2}}{\pi}=\frac{(\frac{1}{4})_{m-n-p}}{(\frac{3}{4})_{n}(\frac{3}{4})_{p}}
\sum_{k=0}^{\infty}(-1)^k
 \frac{(\frac{3}{4})_{k+m}(\frac{3}{4})_{k+n}(\frac{3}{4})_{k+p}}
 {k!(k+m-n)!(k+m-p)!}(3+4m+8k).
 \enm
\end{thm}

Taking $a=\frac{1}{3}+m$, $b=\frac{1}{3}+n$, $c=\frac{1}{3}+p$ in
\eqref{hypergeotric-a} and then letting $d\to\infty$, we gain the
equation.

\begin{thm}\label{thm-ee} For $m,n,p\in \mathbb{Z}$ with $\min\{m-n,m-p,m-2n-2p\}\geq0$,
 there holds the Ramanujan-type series for $1/\pi$ with three free parameters:
\bnm\qquad
 \frac{3\sqrt{3}}{2\pi}=\frac{(\frac{2}{3})_{m-n-p}}{(\frac{1}{3})_{n}(\frac{1}{3})_{p}}
\sum_{k=0}^{\infty}(-1)^k
 \frac{(\frac{1}{3})_{k+m}(\frac{1}{3})_{k+n}(\frac{1}{3})_{k+p}}
 {k!(k+m-n)!(k+m-p)!}(1+3m+6k).
 \enm
\end{thm}

Choosing $a=\frac{2}{3}+m$, $b=\frac{2}{3}+n$, $c=\frac{2}{3}+p$ in
\eqref{hypergeotric-a} and then letting $d\to\infty$, we achieve the
equation.

\begin{thm}\label{thm-ff} For $m,n,p\in \mathbb{Z}$ with $\min\{m-n,m-p,m-2n-2p-1\}\geq0$,
 there holds the Ramanujan-type series for $1/\pi$ with three free parameters:
\bnm\qquad
 \frac{3\sqrt{3}}{2\pi}=\frac{(\frac{1}{3})_{m-n-p}}{(\frac{2}{3})_{n}(\frac{2}{3})_{p}}
\sum_{k=0}^{\infty}(-1)^k
 \frac{(\frac{2}{3})_{k+m}(\frac{2}{3})_{k+n}(\frac{2}{3})_{k+p}}
 {k!(k+m-n)!(k+m-p)!}(2+3m+6k).
 \enm
\end{thm}

Making $a=\frac{1}{12}+m$, $b=\frac{1}{12}+n$, $c=\frac{1}{12}+p$ in
\eqref{hypergeotric-a} and then letting $d\to\infty$, we attain the
equation.

\begin{thm}\label{thm-gg} For $m,n,p\in \mathbb{Z}$ with $\min\{m-n,m-p,m-2n-2p+1\}\geq0$,
 there holds the Ramanujan-type series for $1/\pi$ with three free parameters:
\bnm\qquad
 \frac{3(\sqrt{6}-\sqrt{2})}{\pi}=\frac{(\frac{11}{12})_{m-n-p}}{(\frac{1}{12})_{n}(\frac{1}{12})_{p}}
\sum_{k=0}^{\infty}(-1)^k
 \frac{(\frac{1}{12})_{k+m}(\frac{1}{12})_{k+n}(\frac{1}{12})_{k+p}}
 {k!(k+m-n)!(k+m-p)!}(1+12m+24k).
 \enm
\end{thm}

Setting $a=\frac{5}{12}+m$, $b=\frac{5}{12}+n$, $c=\frac{5}{12}+p$
in \eqref{hypergeotric-a} and then letting $d\to\infty$, we get the
equation.

\begin{thm}\label{thm-hh} For $m,n,p\in \mathbb{Z}$ with $\min\{m-n,m-p,m-2n-2p\}\geq0$,
 there holds the Ramanujan-type series for $1/\pi$ with three free parameters:
\bnm\qquad
 \frac{3(\sqrt{6}+\sqrt{2})}{\pi}=\frac{(\frac{7}{12})_{m-n-p}}{(\frac{5}{12})_{n}(\frac{5}{12})_{p}}
\sum_{k=0}^{\infty}(-1)^k
 \frac{(\frac{5}{12})_{k+m}(\frac{5}{12})_{k+n}(\frac{5}{12})_{k+p}}
 {k!(k+m-n)!(k+m-p)!}(5+12m+24k).
 \enm
\end{thm}

Taking $a=\frac{7}{12}+m$, $b=\frac{7}{12}+n$, $c=\frac{7}{12}+p$ in
\eqref{hypergeotric-a} and then letting $d\to\infty$, we gain the
equation.

\begin{thm}\label{thm-ii} For $m,n,p\in \mathbb{Z}$ with $\min\{m-n,m-p,m-2n-2p\}\geq0$,
 there holds the Ramanujan-type series for $1/\pi$ with three free parameters:
\bnm\qquad
 \frac{3(\sqrt{6}+\sqrt{2})}{\pi}=\frac{(\frac{5}{12})_{m-n-p}}{(\frac{7}{12})_{n}(\frac{7}{12})_{p}}
\sum_{k=0}^{\infty}(-1)^k
 \frac{(\frac{7}{12})_{k+m}(\frac{7}{12})_{k+n}(\frac{7}{12})_{k+p}}
 {k!(k+m-n)!(k+m-p)!}(7+12m+24k).
 \enm
\end{thm}

Choosing $a=\frac{11}{12}+m$, $b=\frac{11}{12}+n$,
$c=\frac{11}{12}+p$ in \eqref{hypergeotric-a} and then letting
$d\to\infty$, we achieve the equation.

\begin{thm}\label{thm-jj} For $m,n,p\in \mathbb{Z}$ with $\min\{m-n,m-p,m-2n-2p-1\}\geq0$,
 there holds the Ramanujan-type series for $1/\pi$ with three free parameters:
\bnm\qquad
 \frac{3(\sqrt{6}-\sqrt{2})}{\pi}=\frac{(\frac{1}{12})_{m-n-p}}{(\frac{11}{12})_{n}(\frac{11}{12})_{p}}
\sum_{k=0}^{\infty}(-1)^k
 \frac{(\frac{11}{12})_{k+m}(\frac{11}{12})_{k+n}(\frac{11}{12})_{k+p}}
 {k!(k+m-n)!(k+m-p)!}(11+12m+24k).
 \enm
\end{thm}

Making $a=\frac{1}{4}+m$, $b=\frac{1}{4}+n$, $c=\frac{1}{4}+p$ and
$d=\frac{1}{2}+q$ in \eqref{hypergeotric-a}, we attain the
 equation.

\begin{thm}\label{thm-a} For $m,n,p,q\in \mathbb{Z}$ with $\min\{m-n,m-p,m-n-p-q\}\geq0$,
 there holds the Ramanujan-type series for $1/\pi$ with four free parameters:
\bnm \quad
 \frac{4}{\pi}&&\xqdn=\:\frac{(\frac{3}{4})_{m-n-p}(\frac{1}{2})_{m-n-q}(\frac{1}{2})_{m-p-q}}
 {(\frac{1}{4})_{m-n-p-q}(\frac{1}{4})_{n}(\frac{1}{4})_{p}(\frac{1}{2})_{q}}
\\&&\xqdn\times\:\: \sum_{k=0}^{\infty}
 \frac{(\frac{1}{4})_{k+m}(\frac{1}{4})_{k+n}(\frac{1}{4})_{k+p}(\frac{1}{2})_{k+q}}
 {k!(k+m-n)!(k+m-p)!(\frac{3}{4})_{k+m-q}}(1+4m+8k).
 \enm
\end{thm}

Setting $a=\frac{3}{4}+m$, $b=\frac{3}{4}+n$, $c=\frac{3}{4}+p$ and
$d=\frac{1}{2}+q$ in \eqref{hypergeotric-a}, we get the equation.

\begin{thm}\label{thm-b} For $m,n,p,q\in \mathbb{Z}$ with $\min\{m-n,m-p,m-n-p-q-1\}\geq0$,
 there holds the Ramanujan-type series for $1/\pi$ with four free parameters:
\bnm \quad
 \frac{1}{\pi}&&\xqdn=\:\frac{(\frac{1}{4})_{m-n-p}(\frac{1}{2})_{m-n-q}(\frac{1}{2})_{m-p-q}}
 {(\frac{3}{4})_{m-n-p-q-1}(\frac{3}{4})_{n}(\frac{3}{4})_{p}(\frac{1}{2})_{q}}
 \\&&\xqdn\times\:\:
\sum_{k=0}^{\infty}
 \frac{(\frac{3}{4})_{k+m}(\frac{3}{4})_{k+n}(\frac{3}{4})_{k+p}(\frac{1}{2})_{k+q}}
 {k!(k+m-n)!(k+m-p)!(\frac{5}{4})_{k+m-q}}(3+4m+8k).
 \enm
\end{thm}

Taking $a=\frac{1}{2}+m$, $b=\frac{1}{2}+n$, $c=\frac{1}{4}+p$ and
$d=\frac{3}{4}+q$ in \eqref{hypergeotric-a}, we gain the equation.

\begin{thm}\label{thm-c} For $m,n,p,q\in \mathbb{Z}$ with $\min\{m-n,m-n-p-q-1\}\geq0$,
 there holds the Ramanujan-type series for $1/\pi$ with four free parameters:
\bnm
 \frac{1}{2\pi}&&\xqdn=\:\frac{(\frac{3}{4})_{m-n-p}(\frac{1}{4})_{m-n-q}(\frac{1}{2})_{m-p-q}}
 {(m-n-p-q-1)!(\frac{1}{2})_{n}(\frac{1}{4})_{p}(\frac{3}{4})_{q}}
\\&&\xqdn\times\:\: \sum_{k=0}^{\infty}
 \frac{(\frac{1}{2})_{k+m}(\frac{1}{2})_{k+n}(\frac{1}{4})_{k+p}(\frac{3}{4})_{k+q}}
 {k!(k+m-n)!(\frac{5}{4})_{k+m-p}(\frac{3}{4})_{k+m-q}}(1+2m+4k).
 \enm
\end{thm}

Choosing $a=\frac{1}{3}+m$, $b=\frac{1}{3}+n$, $c=\frac{1}{3}+p$ and
$d=\frac{1}{2}+q$ in \eqref{hypergeotric-a}, we achieve the
equation.

\begin{thm}\label{thm-d} For $m,n,p,q\in \mathbb{Z}$ with $\min\{m-n,m-p,m-n-p-q\}\geq0$,
 there holds the Ramanujan-type series for $1/\pi$ with four free parameters:
\bnm \quad
 \frac{3\sqrt{3}}{\pi}&&\xqdn=\:\frac{(\frac{2}{3})_{m-n-p}(\frac{1}{2})_{m-n-q}(\frac{1}{2})_{m-p-q}}
 {(\frac{1}{6})_{m-n-p-q}(\frac{1}{3})_{n}(\frac{1}{3})_{p}(\frac{1}{2})_{q}}
\\&&\xqdn\times\:\: \sum_{k=0}^{\infty}
 \frac{(\frac{1}{3})_{k+m}(\frac{1}{3})_{k+n}(\frac{1}{3})_{k+p}(\frac{1}{2})_{k+q}}
 {k!(k+m-n)!(k+m-p)!(\frac{5}{6})_{k+m-q}}(1+3m+6k).
 \enm
\end{thm}

Making $a=\frac{2}{3}+m$, $b=\frac{2}{3}+n$, $c=\frac{2}{3}+p$ and
$d=\frac{1}{2}+q$ in \eqref{hypergeotric-a}, we attain the equation.

\begin{thm}\label{thm-e} For $m,n,p,q\in \mathbb{Z}$ with $\min\{m-n,m-p,m-n-p-q-1\}\geq0$,
 there holds the Ramanujan-type series for $1/\pi$ with four free parameters:
\bnm \qquad
 \frac{\sqrt{3}}{2\pi}&&\xqdn=\:\frac{(\frac{1}{3})_{m-n-p}(\frac{1}{2})_{m-n-q}(\frac{1}{2})_{m-p-q}}
 {(\frac{5}{6})_{m-n-p-q-1}(\frac{2}{3})_{n}(\frac{2}{3})_{p}(\frac{1}{2})_{q}}
 \\&&\xqdn\times\:\:
\sum_{k=0}^{\infty}
 \frac{(\frac{2}{3})_{k+m}(\frac{2}{3})_{k+n}(\frac{2}{3})_{k+p}(\frac{1}{2})_{k+q}}
 {k!(k+m-n)!(k+m-p)!(\frac{7}{6})_{k+m-q}}(2+3m+6k).
 \enm
\end{thm}

Setting $a=\frac{1}{6}+m$, $b=\frac{1}{6}+n$, $c=\frac{1}{6}+p$ and
$d=\frac{1}{2}+q$ in \eqref{hypergeotric-a}, we get the equation.

\begin{thm}\label{thm-g} For $m,n,p,q\in \mathbb{Z}$ with $\min\{m-n,m-p,m-n-p-q\}\geq0$,
 there holds the Ramanujan-type series for $1/\pi$ with four free parameters:
\bnm \quad
 \frac{2\sqrt{3}}{\pi}&&\xqdn=\:\frac{(\frac{5}{6})_{m-n-p}(\frac{1}{2})_{m-n-q}(\frac{1}{2})_{m-p-q}}
 {(\frac{1}{3})_{m-n-p-q}(\frac{1}{6})_{n}(\frac{1}{6})_{p}(\frac{1}{2})_{q}}
\\&&\xqdn\times\:\: \sum_{k=0}^{\infty}
 \frac{(\frac{1}{6})_{k+m}(\frac{1}{6})_{k+n}(\frac{1}{6})_{k+p}(\frac{1}{2})_{k+q}}
 {k!(k+m-n)!(k+m-p)!(\frac{2}{3})_{k+m-q}}(1+6m+12k).
 \enm
\end{thm}

Taking $a=\frac{5}{6}+m$, $b=\frac{5}{6}+n$, $c=\frac{5}{6}+p$ and
$d=\frac{1}{2}+q$ in \eqref{hypergeotric-a}, we gain the equation.

\begin{thm}\label{thm-h} For $m,n,p,q\in \mathbb{Z}$ with $\min\{m-n,m-p,m-n-p-q-1\}\geq0$,
 there holds the Ramanujan-type series for $1/\pi$ with four free parameters:
\bnm \quad
 \frac{2\sqrt{3}}{3\pi}&&\xqdn=\:\frac{(\frac{1}{6})_{m-n-p}(\frac{1}{2})_{m-n-q}(\frac{1}{2})_{m-p-q}}
 {(\frac{2}{3})_{m-n-p-q-1}(\frac{5}{6})_{n}(\frac{5}{6})_{p}(\frac{1}{2})_{q}}
 \\&&\xqdn\times\:\:
\sum_{k=0}^{\infty}
 \frac{(\frac{5}{6})_{k+m}(\frac{5}{6})_{k+n}(\frac{5}{6})_{k+p}(\frac{1}{2})_{k+q}}
 {k!(k+m-n)!(k+m-p)!(\frac{4}{3})_{k+m-q}}(5+6m+12k).
 \enm
\end{thm}

Choosing $a=\frac{1}{2}+m$, $b=\frac{1}{2}+n$, $c=\frac{1}{6}+p$ and
$d=\frac{5}{6}+q$ in \eqref{hypergeotric-a}, we achieve the
equation.

\begin{thm}\label{thm-i} For $m,n,p,q\in \mathbb{Z}$ with $\min\{m-n,m-n-p-q-1\}\geq0$,
 there holds the Ramanujan-type series for $1/\pi$ with four free parameters:
\bnm
 \frac{2\sqrt{3}}{9\pi}&&\xqdn=\:\frac{(\frac{5}{6})_{m-n-p}(\frac{1}{6})_{m-n-q}(\frac{1}{2})_{m-p-q}}
 {(m-n-p-q-1)!(\frac{1}{2})_{n}(\frac{1}{6})_{p}(\frac{5}{6})_{q}}
\\&&\xqdn\times\:\: \sum_{k=0}^{\infty}
 \frac{(\frac{1}{2})_{k+m}(\frac{1}{2})_{k+n}(\frac{1}{6})_{k+p}(\frac{5}{6})_{k+q}}
 {k!(k+m-n)!(\frac{4}{3})_{k+m-p}(\frac{2}{3})_{k+m-q}}(1+2m+4k).
 \enm
\end{thm}

 Making $a=\frac{1}{12}+m$, $b=\frac{1}{12}+n$, $c=\frac{1}{12}+p$
 and $d=\frac{1}{2}+q$ in \eqref{hypergeotric-a}, we attain the
 equation.

\begin{thm}\label{thm-j} For $m,n,p,q\in \mathbb{Z}$ with $\min\{m-n,m-p,m-n-p-q\}\geq0$,
 there holds the Ramanujan-type series for $1/\pi$ with four free parameters:
\bnm
 \frac{12(2-\sqrt{3})}{\pi}&&\xqdn=\:\frac{(\frac{11}{12})_{m-n-p}(\frac{1}{2})_{m-n-q}(\frac{1}{2})_{m-p-q}}
 {(\frac{5}{12})_{m-n-p-q}(\frac{1}{12})_{n}(\frac{1}{12})_{p}(\frac{1}{2})_{q}}
\\&&\xqdn\times\:\: \sum_{k=0}^{\infty}
 \frac{(\frac{1}{12})_{k+m}(\frac{1}{12})_{k+n}(\frac{1}{12})_{k+p}(\frac{1}{2})_{k+q}}
 {k!(k+m-n)!(k+m-p)!(\frac{7}{12})_{k+m-q}}(1+12m+24k).
 \enm
\end{thm}

Setting $a=\frac{5}{12}+m$, $b=\frac{5}{12}+n$, $c=\frac{5}{12}+p$
and $d=\frac{1}{2}+q$ in \eqref{hypergeotric-a}, we get the
 equation.

\begin{thm}\label{thm-k} For $m,n,p,q\in \mathbb{Z}$ with $\min\{m-n,m-p,m-n-p-q\}\geq0$,
 there holds the Ramanujan-type series for $1/\pi$ with four free parameters:
\bnm
 \frac{12(2+\sqrt{3})}{\pi}&&\xqdn=\:\frac{(\frac{7}{12})_{m-n-p}(\frac{1}{2})_{m-n-q}(\frac{1}{2})_{m-p-q}}
 {(\frac{1}{12})_{m-n-p-q}(\frac{5}{12})_{n}(\frac{5}{12})_{p}(\frac{1}{2})_{q}}
\\&&\xqdn\times\:\: \sum_{k=0}^{\infty}
 \frac{(\frac{5}{12})_{k+m}(\frac{5}{12})_{k+n}(\frac{5}{12})_{k+p}(\frac{1}{2})_{k+q}}
 {k!(k+m-n)!(k+m-p)!(\frac{11}{12})_{k+m-q}}(5+12m+24k).
 \enm
\end{thm}

Taking $a=\frac{7}{12}+m$, $b=\frac{7}{12}+n$, $c=\frac{7}{12}+p$
and $d=\frac{1}{2}+q$ in \eqref{hypergeotric-a}, we gain the
equation.

\begin{thm}\label{thm-l} For $m,n,p,q\in \mathbb{Z}$ with $\min\{m-n,m-p,m-n-p-q-1\}\geq0$,
 there holds the Ramanujan-type series for $1/\pi$ with four free parameters:
\bnm
 \frac{2+\sqrt{3}}{\pi}&&\xqdn=\:\frac{(\frac{5}{12})_{m-n-p}(\frac{1}{2})_{m-n-q}(\frac{1}{2})_{m-p-q}}
 {(\frac{11}{12})_{m-n-p-q-1}(\frac{7}{12})_{n}(\frac{7}{12})_{p}(\frac{1}{2})_{q}}
\\&&\xqdn\times\:\: \sum_{k=0}^{\infty}
 \frac{(\frac{7}{12})_{k+m}(\frac{7}{12})_{k+n}(\frac{7}{12})_{k+p}(\frac{1}{2})_{k+q}}
 {k!(k+m-n)!(k+m-p)!(\frac{13}{12})_{k+m-q}}(7+12m+24k).
 \enm
\end{thm}

Choosing $a=\frac{11}{12}+m$, $b=\frac{11}{12}+n$,
$c=\frac{11}{12}+p$ and $d=\frac{1}{2}+q$ in \eqref{hypergeotric-a},
we achieve the equation.

\begin{thm}\label{thm-m} For $m,n,p,q\in \mathbb{Z}$ with $\min\{m-n,m-p,m-n-p-q-1\}\geq0$,
 there holds the Ramanujan-type series for $1/\pi$ with four free parameters:
\bnm
 \frac{5(2-\sqrt{3})}{\pi}&&\xqdn=\:\frac{(\frac{1}{12})_{m-n-p}(\frac{1}{2})_{m-n-q}(\frac{1}{2})_{m-p-q}}
 {(\frac{7}{12})_{m-n-p-q-1}(\frac{11}{12})_{n}(\frac{11}{12})_{p}(\frac{1}{2})_{q}}
\\&&\xqdn\times\:\: \sum_{k=0}^{\infty}
 \frac{(\frac{11}{12})_{k+m}(\frac{11}{12})_{k+n}(\frac{11}{12})_{k+p}(\frac{1}{2})_{k+q}}
 {k!(k+m-n)!(k+m-p)!(\frac{17}{12})_{k+m-q}}(11+12m+24k).
 \enm
\end{thm}

Making $a=\frac{1}{2}+m$, $b=\frac{1}{2}+n$, $c=\frac{1}{12}+p$ and
$d=\frac{11}{12}+q$ in \eqref{hypergeotric-a}, we attain the
equation.

\begin{thm}\label{thm-n} For $m,n,p,q\in \mathbb{Z}$ with $\min\{m-n,m-n-p-q-1\}\geq0$,
 there holds the Ramanujan-type series for $1/\pi$ with four free parameters:
\bnm
 \frac{5(2-\sqrt{3})}{6\pi}&&\xqdn=\:\frac{(\frac{11}{12})_{m-n-p}(\frac{1}{12})_{m-n-q}(\frac{1}{2})_{m-p-q}}
 {(m-n-p-q-1)!(\frac{1}{2})_{n}(\frac{1}{12})_{p}(\frac{11}{12})_{q}}
\\&&\xqdn\times\:\: \sum_{k=0}^{\infty}
 \frac{(\frac{1}{2})_{k+m}(\frac{1}{2})_{k+n}(\frac{1}{12})_{k+p}(\frac{11}{12})_{k+q}}
 {k!(k+m-n)!(\frac{17}{12})_{k+m-p}(\frac{7}{12})_{k+m-q}}(1+2m+4k).
 \enm
\end{thm}

Setting $a=\frac{1}{2}+m$, $b=\frac{1}{2}+n$, $c=\frac{5}{12}+p$ and
$d=\frac{7}{12}+q$ in \eqref{hypergeotric-a}, we get the equation.

\begin{thm}\label{thm-o} For $m,n,p,q\in \mathbb{Z}$ with $\min\{m-n,m-n-p-q-1\}\geq0$,
 there holds the Ramanujan-type series for $1/\pi$ with four free parameters:
\bnm
 \frac{2+\sqrt{3}}{6\pi}&&\xqdn=\:\frac{(\frac{7}{12})_{m-n-p}(\frac{5}{12})_{m-n-q}(\frac{1}{2})_{m-p-q}}
 {(m-n-p-q-1)!(\frac{1}{2})_{n}(\frac{5}{12})_{p}(\frac{7}{12})_{q}}
\\&&\xqdn\times\:\: \sum_{k=0}^{\infty}
 \frac{(\frac{1}{2})_{k+m}(\frac{1}{2})_{k+n}(\frac{5}{12})_{k+p}(\frac{7}{12})_{k+q}}
 {k!(k+m-n)!(\frac{13}{12})_{k+m-p}(\frac{11}{12})_{k+m-q}}(1+2m+4k).
 \enm
\end{thm}

 \textbf{Remark:} With the change of the parameters, Theorems
\ref{thm-aa}-\ref{thm-o} can produce numerous concrete
 Ramanujan-type series for $1/\pi$. We shall not lay them out one by one here.
%%%%%%%%%%%%%%%%%%%%%%%%%%%%%%%%%%%%%%%%%%%%%%%%%%%%%%%%%%%%%%%%%%%

%%%%%%%%%%%%%%%%%%%%%%%%%%%%%%%%%%%%%%%%%%%%%%%%%%%%%%%%%%%%%%%%%%%
%%%%%%%%%%%%%%%%%%%%%%%%%%%%%%%%%%%%%%%%%%%%%%%%%%%%%%%%%%%%%%%%%%%
%%%%%%%%%%%%%%%%%%%%%%%%%%%%%%%%%%%%%%%%%%%%%%%%%%%%%%%%%%%%%%%%%%%

\end{document}